\documentclass[twoside,leqno,twocolumn]{article}

\usepackage[letterpaper]{geometry}

\usepackage{ltexpprt}
\usepackage{multirow, multicol}
\usepackage{graphicx}
\usepackage{caption}
\usepackage{rotating}
\usepackage{booktabs} 
\usepackage{comment}
\usepackage{amsfonts}  
\usepackage{hyperref}
\usepackage{algorithm, algpseudocode}
\usepackage{adjustbox}
\usepackage{subcaption}
\algnewcommand{\algorithmicforeach}{\textbf{for each}}
\algdef{SE}[FOR]{ForEach}{EndForEach}[1]
  {\algorithmicforeach\ #1\ \algorithmicdo}
  {\algorithmicend\ \algorithmicforeach}
\usepackage{amsmath}
  
\usepackage{cleveref}
\usepackage{booktabs}

\usepackage{xcolor}

\newcommand{\tofill}[1]{{\color{red} [FILL]}}

\def\Z{\mathbb{Z}}
\def\R{\mathbb{R}}
\def\C{\mathbb{C}}
\def\id{\mathrm{Id}}
\def\1{\mathbf{1}}
\newcommand{\dualdot}[2]{\left\langle#1,#2\right\rangle}

\def\qupid{\textsc{Qupid}}
\def\Id{\texttt{id}}
\def\fft{\texttt{fft}}
\newcommand{\db}[1]{\texttt{db#1}}
\newcommand{\coif}[1]{\texttt{coif#1}}
\def\orbit5k{\texttt{ORBIT5K}}
\def\CDs{\texttt{CD68}$^+$}
\def\CDh{\texttt{CD8}$^+$}
\def\Fox{\texttt{FoxP3}$^+$}

\DeclareMathOperator{\mat}{\mathcal{M}}
\DeclareMathOperator{\meas}{\mathcal{M}}
\DeclareMathOperator{\Fourier}{\mathcal{F}}
\def\Diag{\mathcal{D}}
\def\Quant{\mathsf{q}}
\newcommand{\quant}[1]{\mathsf{q}\mu_{#1}}
\DeclareMathOperator{\HKS}{HKS}

\begin{document}
\onecolumn

\newcommand\relatedversion{}

\title{\Large Discrete transforms of quantized persistence diagrams\relatedversion}

 \author{
    Michael E. Van Huffel\thanks{Department of Mathematics, ETH Zurich, Switzerland, \texttt{michavan@student.ethz.ch}}
    \and
    Olympio Hacquard\thanks{KUIAS, Kyoto University, Japan, \texttt{hacquard.olympio.47i@st.kyoto-u.ac.jp}}
     \and
    Vadim Lebovici\thanks{Mathematical Institute, University of Oxford, United Kingdom, \texttt{vadim.lebovici@maths.ox.ac.uk}}
    \and Matteo Palo\thanks{Department of Mathematics, ETH Zurich, Switzerland, \texttt{mapalo@student.ethz.ch}}
 }

\date{}
\twocolumn
\maketitle


\fancyfoot[R]{\scriptsize{Copyright \textcopyright\ 2025 by SIAM\\
Unauthorized reproduction of this article is prohibited}}





\begin{abstract}
    Topological data analysis leverages topological features to analyze datasets, with applications in diverse fields like medical sciences and biology. A key tool of this theory is the persistence diagram, which encodes topological information but poses challenges for integration into standard machine learning pipelines. 
    We introduce \qupid{} (QUantized Persistence and Integral transforms of Diagrams), a novel and simple method for vectorizing persistence diagrams. First, \qupid{} uses a binning procedure to turn persistence diagrams into finite measures on a grid and then applies discrete transforms to these measures. Key features are the choice of log-scaled grids that emphasize information contained near the diagonal in persistence diagrams, combined with the use of discrete transforms to enhance and efficiently encode the obtained topological information.
    We conduct an in-depth experimental analysis of \qupid{}, showing that the simplicity of our method results in very low computational costs while preserving highly competitive performances compared to state-of-the-art methods across numerous classification tasks on both synthetic and real-world datasets. 
    Finally, we provide experimental evidence that our method is robust to a decrease in the grid resolution used.
\end{abstract}

The code used to reproduce all the experiments is available at \href{https://github.com/majkevh/qupid}{\texttt{https://github.com/majkevh/qupid}}.

\section{Introduction} \label{sec:intro}

Topological data analysis is a branch of machine learning that consists of extracting topological and geometric information to help analyze datasets. This information can be easily included in a workflow to perform various learning tasks on the
data. As a result, these techniques have found applications in fields as diverse as medical sciences \cite{rieck2020uncovering, fernandez2022topological, aukerman2022persistent}, biology, \cite{ichinomiya2020protein, rabadan2019topological}, material sciences, \cite{lee2017quantifying, hiraoka2016hierarchical}, cosmology \cite{pranav2017topology}, music theory \cite{alcala2024framework, mijangos2024musical} and many more. A central topological descriptor of this theory is the persistence diagram. Mathematically, a persistence diagram is a point cloud---that is, a finite subset---of $\R^2$, where each point relates to and gives quantitative information on a topological feature of the data. Persistence diagrams can be computed from a wide spectrum of input data, including point clouds, graphs, images or time series. One natural approach is to replace raw input data with their persistence diagrams in classical machine learning analyses---such as regression or classification tasks---to focus on the topological and geometric information contained in data.

Due to their point cloud nature, persistence diagrams cannot be straightforwardly embedded in most standard machine learning algorithms. Several methods have been proposed to overcome this difficulty, such as neural network architectures pioneered in \cite{qi2017pointnet}, or works based on optimal transport as in \cite{moosmuller2023linear,khurana2023supervised}. 

Due to the topological information they carry, there has been a growing interest in proposing point cloud learning methods specifically tailored to persistence diagrams. Apart from some methods that perform learning directly at the point cloud level such as \cite{bba}, most aim at mapping persistence diagrams into a vector space. Following this vectorization step, any off-the-shelf machine learning algorithm suited to the task to perform (e.g., classification, regression, or clustering) can be used. 
In \cite{adams2017persistence} for instance, persistence diagrams are turned into so-called \textit{persistence images}, which consist of the pixelization of a convolution of the point cloud with a Gaussian kernel. Some methods such as \cite{carriere2017sliced} and \cite{le2018persistence} have used kernel-based algorithms. More recently, neural network-based vectorizations have been proposed, such as \cite{pmlr-v108-carriere20a} and \cite{reinauer2021persformer}. Finally, let us mention ATOL \cite{pmlr-v130-royer21a}, which uses a clustering algorithm on persistence diagrams to find points of interest and takes the mass around each cluster centre as coordinates for the vectorization.

\paragraph*{Contributions.} We introduce a new vectorization methodology for persistence diagrams called \qupid{} (QUantized Persistence and Integral transforms of Diagrams). First, we turn persistence diagrams into finite measures on a user-chosen grid via a binning procedure. Second, we apply a discrete integral transform (e.g. Fourier, wavelet) to these measures to extract and enhance valuable structure patterns in persistence diagrams. 

We claim that the strength of our method lies in its great simplicity. We show that this simplicity translates into reduced computational time. We obtain very competitive accuracy results in supervised classification tasks while being several times faster than persistence images and kernel methods; and being dramatically faster than neural networks. In particular, we validated our model on a real-world dataset of tumour immune cell type classification. Finally, combining the classical knowledge of discrete transforms and the structure of persistence diagrams, the simplicity of our method allows for straightforward interpretations of the output vector representation.

We interpret the promising results of \qupid{} as due to three key features. First, our binning procedure greatly emphasises regions of persistence diagrams with a high density of points, which typically correspond to regions where fine geometric information on data is experimentally proven to be contained; see \Cref{sec:quantization} for more details. 
Second, we leverage the flexibility of our method in the possibility of choosing log-scaled grids that further emphasise the high-density regions in persistence diagrams. Doing so, \qupid{} becomes a vectorization of point clouds specifically adapted, but not restricted to, persistence diagrams. 
Third, in comparison with the quite similar ATOL algorithm, the possibility to choose any discrete transform allows for additional flexibility and a clear improvement in feature extraction from persistence diagrams. 

\paragraph*{Outline.}
\Cref{sec:background} introduces basic notions of topological data analysis and persistence diagrams, emphasizing why they are relevant data descriptors. \Cref{sec:methodology} describes the algorithm \qupid{} to turn persistence diagrams into vectors. Finally, \Cref{sec:experiments} presents numerous experiments on synthetic and real-world datasets, showing that \qupid{} reaches a very competitive accuracy when benchmarked against state-of-the-art methods. In addition, this section provides a thorough analysis of our method, its computational efficiency, and the influence of the resolution of the grid used in the quantization.

\section{Background on topological data analysis}
\label{sec:background}
Throughout this paper, we develop an algorithm for fast learning on point cloud data. The main emphasis will be the classification of topological descriptors called persistence diagrams. In this section, we describe their construction and build intuition on the type of information they capture.

Starting from diverse input data (graphs, images, point clouds, time series...), topological data analysis aims at extracting topological information, that is, the homotopy type of some underlying geometric object the data is assumed to lie on. The main tool to do so is \emph{homology} from algebraic topology. Doing so with a multi-scale approach is the purpose of \textit{persistent homology}, as introduced in the reference textbooks~\cite{edelsbrunner2022computational,boissonnat2018geometric}.

Intuitively, given a topological space~$X$, the~$k$-th homology of~$X$ is a vector space whose dimension is equal to the number of independent k-dimensional holes of X for~$k \geq 1$ (loops for~$k = 1$, voids for~$k=2$ and so on) and to the number of connected components for~$k=0$.
Given a nested family of topological spaces~$\mathcal{X} = (X_t)_{t \in \mathbb{R}}$ called a \textit{filtration} (i.e., such that~$X_{t_1} \subset X_{t_2}$ whenever~$t_1 < t_2$), the degree-$k$ persistent homology tracks the evolution of the degree-$k$ homological features---i.e., the generators of the~$k$-th homology of~$X_t$ as the parameter~$t$ increases. 

We denote by $\Diag$ the set of finite multi-sets of points~$(a,b)\in\R^2$ such that $a\leq b$, called \emph{persistence diagrams}. To a filtration $\mathcal{X}$ and an integer~$k\geq 0$ one can associate a persistence diagram $D\in\Diag$ composed of the collection of all couples $(a, b)$ such that a given degree-$k$ homological feature is born at parameter~$t=a$ and dies at parameter~$t=b$ in the filtration~$\mathcal{X}$.

Given raw input data, the choice of filtration is central and will determine the aspect of the persistence diagram. If the input data is a point cloud $\mathbf{X} = (x_i)_{i=1}^n$ in $\mathbb{R}^d$, a widespread choice of filtration is given by the union of balls centred around each data point for increasing radius $t\in\R_{\geq 0}$, called the \emph{\v{C}ech filtration} and denoted by $\check{C} (\mathbf{X}, t) = \cup_{i=1}^n B(x_i, t) $ for each $t \in \R_{\geq 0}$. We therefore have $d$ persistence diagrams corresponding to topological features of dimensions $0, 1, \ldots, d-1$. Unless specified otherwise, this is the filtration we will use when processing point cloud data in \Cref{sec:experiments}. 

\begin{figure*}[]
    \centering  
    \begin{subfigure}[]{0.22\linewidth}
        \centering
        \includegraphics[height=3.5cm]{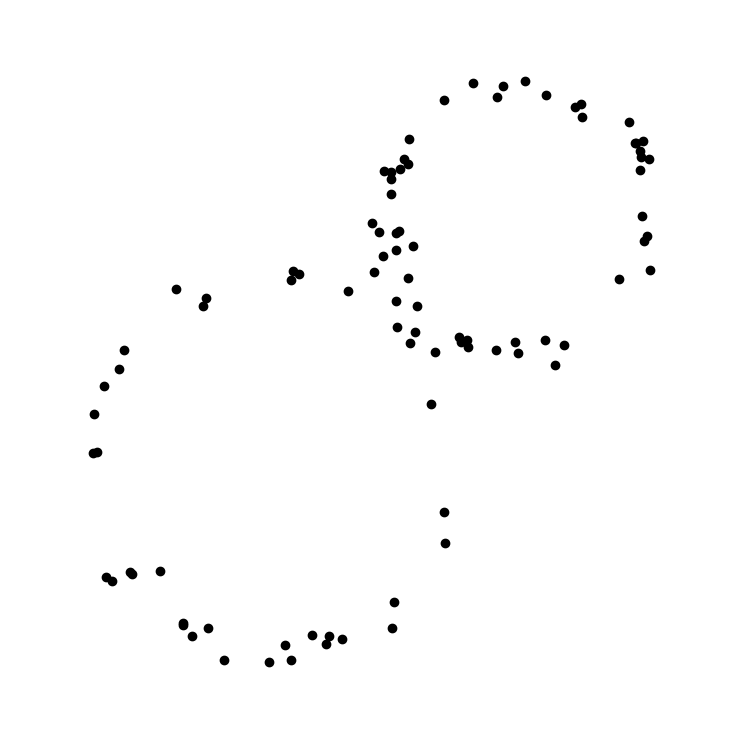}
        \caption{Point cloud $\mathbf{X}$}
        \label{fig:raw_PC}
    \end{subfigure}
    \begin{subfigure}[]{0.22\linewidth}
        \centering
        \includegraphics[height=3.5cm]{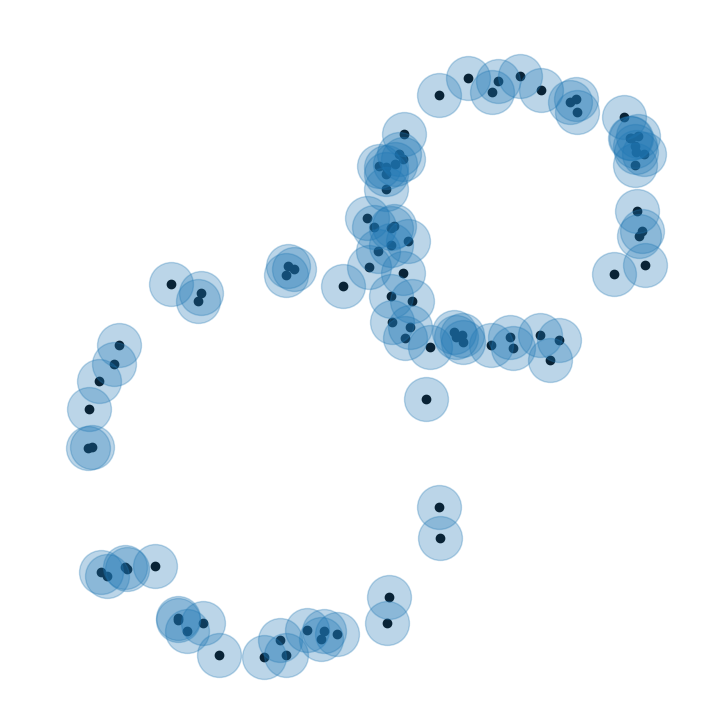}
        \caption{$\check{C} (\mathbf{X}, t_1)$}
        \label{fig:cech_1}
    \end{subfigure}
        \begin{subfigure}[]{0.22\linewidth}
        \centering
        \includegraphics[height=3.5cm]{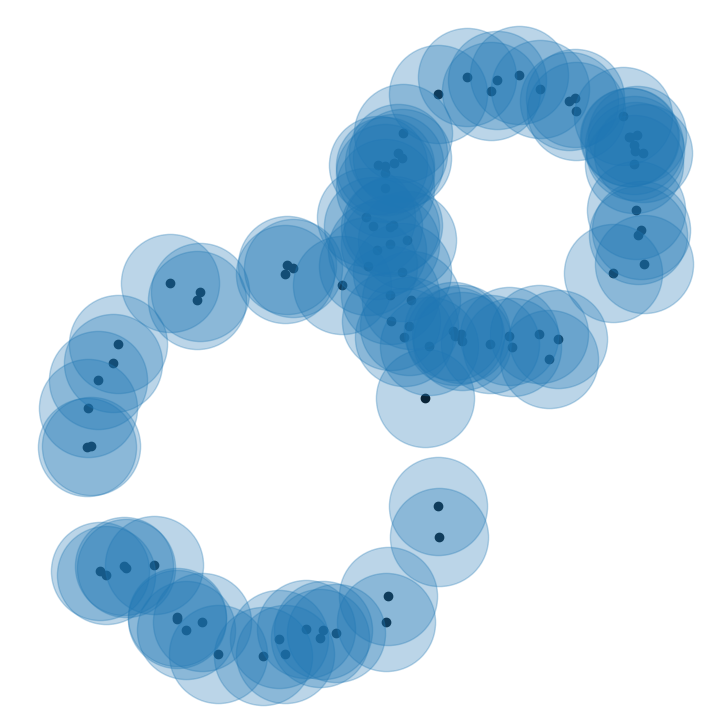}
        \caption{$\check{C} (\mathbf{X}, t_2)$}
        \label{fig:cech_1}
    \end{subfigure}
    \begin{subfigure}[]{0.3\linewidth}
        \centering
        \includegraphics[height=4cm]{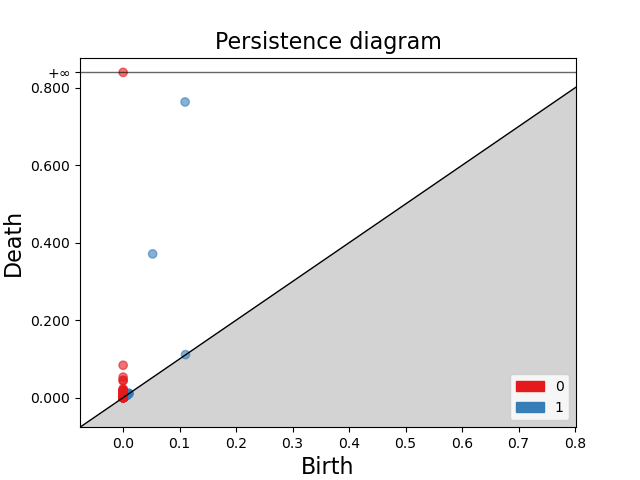}
        \caption{degree-0 (red) and degree-1 (blue) persistence diagrams}
        \label{fig:toy_PD}
        \end{subfigure}
    \caption{Union of balls at radii $t_1$ (b) and $t_2$ (c) with $t_1<t_2$ centred on a point cloud sampled around two circles (a). Taking every radius $t$ defines the \v{C}ech filtration and the corresponding persistence diagrams (d).}
    \label{fig:ex_PD}
\end{figure*}

We show an example of the construction of a persistence diagram for the \v{C}ech filtration of a point cloud of $\mathbb{R}^2$ in \Cref{fig:ex_PD}. In this example, the points are sampled on two circles of different radius. Once the radius in the \v{C}ech filtration is large enough, two cycles appear successively. These topological features exist until the radius parameter becomes large enough and the union of balls fills each circle. Hence, we observe two points significantly away from the diagonal in the degree-$1$ persistence diagram, corresponding to the two circles in the original data, as well as numerous points lying close to the diagonal. We remark that the $0$-persistence diagram is a one-dimensional object. Indeed, there are as many connected components as points in the initial data. They are all born at time~0 and die successively by merging with neighbouring points, such that they can be encoded using only their death time.

We define the \textit{persistence} of a feature $(b, d)$ in the diagram as the real number $d-b$, corresponding to its lifetime. The situation from Figure \ref{fig:ex_PD} is typical. Indeed, it is now well-understood that the persistence diagram of a \v{C}ech filtration is composed of a few points of high persistence on the one hand, corresponding to global homological features of the input point cloud, and a large clump of points of low persistence on the other hand, corresponding to local and sampling effects. We refer to \cite{bba} for a more detailed overview.

Finally, note that as observed in \cite{chazal2013structure}, a persistence diagram~$D\in\Diag$ can be equivalently expressed as a discrete measure~$\mu_D$ over~$\R^2$, called \emph{persistence measure} of~$D$, defined as the sum of Dirac masses centred on each point of the diagram~$\mu_D = \sum_{(b,d) \in D} \delta_{(b, d)}$. This measure-theoretic approach has led to many developments in the structure of the space of persistence diagrams, see \cite{divol2019understanding}. 

\section{Methodology}\label{sec:methodology}

An objective of topological data analysis is to extract predictions from the persistence diagrams. Unfortunately, given their mathematical nature as a multi-set of points (or a discrete measure), they cannot directly replace the raw data as input of standard machine learning algorithms. Indeed, they almost always require the data to be elements of a vector space. A common strategy is therefore to map the persistence diagrams in $\mathbb{R}^n$ where the dimension $n$ is typically chosen by the user. This initial data-transformation step is commonly referred to as the \textit{vectorization step}. We refer to \cite{hensel2021survey} for a survey of different vectorization methods in topological machine learning.

In this section, we introduce the \qupid{} vectorization algorithm. It proceeds in two steps:
\begin{equation*}
    \qupid\colon \Diag \overset{\Quant}{\longrightarrow} \meas_G \overset{\Psi}{\longrightarrow} \C^n,
\end{equation*}
where the first step~$\Quant\colon \Diag \to \mat_G$ is a quantization step sending a persistence diagram to a finite measure supported on a fixed grid~$G$ (\Cref{sec:quantization}), and the second step~$\Psi: \meas_G \to \C^n$ is a discrete transform (\Cref{sec:discrete-transforms}) that outputs a vector in a Euclidean space of dimension $n$. See \Cref{fig:methodo} for an illustration.

\subsection{Quantized persistence diagrams}\label{sec:quantization}
Let $D\in \Diag$ be a persistence diagram. The \emph{persistence diagram in birth-persistence coordinates} associated to $D$ is the image $D^\bot$ of $D$ by the map $(b,d)\mapsto (b,d-b)$. The grid used for quantized persistence diagrams will be associated with two finite sets of non-negative real numbers~$\mathcal{B} = \{b_1 <\ldots<b_r\}$ and $\mathcal{P}=\{p_1 < \ldots < p_s\}$. Denote~$G = \{1,\ldots,r\}\times\{1,\ldots, s\}$ and by~$\meas_G$ the space of finite measures supported on~$G$. The quantization of persistence diagrams is a simple binning process:
\begin{definition}
    The \emph{quantization} of~$D$ with respect to~$(\mathcal{B},\mathcal{P})$ is the finite measure~$\quant{D} \in \meas_G$ defined by:
    \begin{equation*}
        \quant{D} = \sum_{(i,j)\in G} m_{i,j} \delta_{(i,j)},
    \end{equation*}
    where, setting~$b_{r+1} = p_{s+1} = +\infty$ and denoting by~$\1_Z \colon \R^2 \to \{0,1\}$ the indicating function of $Z\subset \R^2$,
    \begin{align*}
        m_{i,j} &= \dualdot{\mu_D}{\1_{[b_i,b_{i+1})\times[p_j,p_{j+1})}} \\
        &= \# \big\{(b,p)\in D^\bot \colon b_{i}\leq b< b_{i+1} \\
        &\hspace{2.5cm}\mbox{ and } p_{j}\leq p < p_{j+1} \big\}.
    \end{align*} 
\end{definition}

\begin{figure*}[t]
    \centering  
    \begin{subfigure}[]{0.24\linewidth}
        \centering
        \includegraphics[width=1\linewidth]{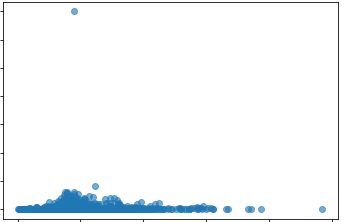}
        \caption{}
        \label{fig:methodo_pd}
    \end{subfigure}
    \begin{subfigure}[]{0.24\linewidth}
        \centering
        \includegraphics[width=\linewidth]{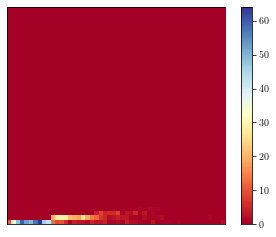}
        \caption{}
        \label{fig:methodo_grid_id}
    \end{subfigure}
        \begin{subfigure}[]{0.24\linewidth}
        \centering
        \includegraphics[width=\linewidth]{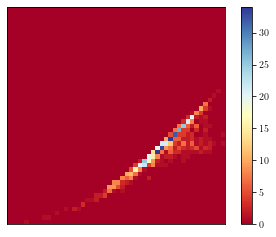}
        \caption{}
        \label{fig:methodo_grid}
    \end{subfigure}
    \begin{subfigure}[]{0.24\linewidth}
        \centering
        \includegraphics[width=\linewidth]{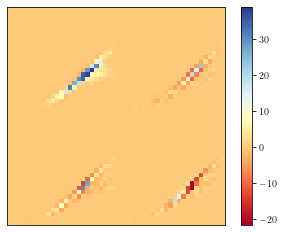}
        \caption{}
    \end{subfigure}
    \caption{An example of degree-1 persistence diagram in birth-persistence coordinates obtained from the \v{C}ech filtration of a point cloud (a), its associated quantized persistence diagram with regular grid of size $50\times 50$ (b), its associated quantized diagram with log-scaled grid of size $50\times 50$ and parameter $\alpha = (500,500)$ (c), and the discrete Daubechies transform of order $2$ (with approximation coefficients, vertical, horizontal, diagonal details from top left to bottom right) of this quantized diagram (d).}
    \label{fig:methodo}
\end{figure*}

In our experiments, we will exploit the flexibility of our methodology concerning the choice of finite sets~$\mathcal{B} = \{b_1 <\ldots<b_r\}$ and $\mathcal{P}=\{p_1 < \ldots < p_s\}$. The first obvious choice is to uniformly sample values between the extremal values attained by the persistence diagrams in birth-persistence coordinates computed from data. We refer to this choice as a \emph{uniform grid}. The second one, called \emph{log-scaled grid}, is to apply a scaling function to the birth-persistence space before the quantization step, namely the function~$(b,p) \mapsto \log(1 + b\alpha_1 + p\alpha_2)$ for some~$\alpha = (\alpha_1, \alpha_2)$ with $\alpha_1, \alpha_2 > 0$. See~\Cref{fig:methodo_pd,fig:methodo_grid_id,fig:methodo_grid} for an illustration. 

Our methodology is designed to emphasize low persistence points. First, as explained in \Cref{sec:background}, persistence diagrams typically contain few high-persistence points and numerous low-persistence points. Thus, as illustrated in \Cref{fig:methodo}, the few high persistence points in the persistence diagrams will be significantly shaded in the binning process by the high multiplicities of low persistence regions. Moreover, the ability to choose an appropriate rescaling parameter~$\alpha$ in log-scaled grids provides quantized persistence diagrams with finer resolution in low persistence regions. 

This emphasis on low-persistence points will be key to our experimental results in \Cref{sec:experiments}. Indeed, if historically considered noise as in \cite{adams2017persistence}, low persistence points are now considered as containing crucial information on data, see \cite{bubenik20} and \cite{bba}. In point cloud analysis, for instance, low persistence points are tightly connected to pairwise distances between points of the point cloud, retaining a lot of geometric and sampling-related information. These considerations have been used in curvature regression tasks \cite{bubenik20,reinauer2021persformer,hacquard23} or in the unsupervised classification of geometric patterns in background noise~\cite{hacquard23}. The experiments in \Cref{sec:experiments} aim to provide new evidence in this direction.

Remark that our method is quite similar to the ATOL method from \cite{pmlr-v130-royer21a}. One of the main differences lies in the way low-persistence features are emphasized: we chose an adaptive-grid approach where the scaling parameter is optimized during the learning process. On the contrary, ATOL runs a clustering algorithm on the data such that their binning process is fully independent of the data's labels. In addition, taking a large class of integral transforms allows us for a more diverse data exploration, that can result in a quite strong accuracy gain in some tasks. We refer in particular to Section \ref{sec:cancer}.

\subsection{Discrete transforms}\label{sec:discrete-transforms} In this section, we briefly recall some classical notions from the theory of discrete transforms and refer to \cite{mallat09} for further details.
Consider a finite grid $G = \{1,\ldots,r\}\times\{1,\ldots,s\}\subset \R^2$ as in \Cref{sec:quantization}. We denote by $n = r \times s$ the resolution of the grid. The grid is chosen by the user and its resolution $n$ has a central importance as it acts as a trade-off between the computational efficiency and numerical accuracy of the method. We let $\{\psi_l\}_{1\leq l\leq n}$ be a family of square-integrable complex-valued functions on $\R^2$. 

\begin{definition}\label{def:transforms}
    The \emph{discrete transform} associated to $\{\psi_l\}$ of the finite measure~$\nu\in \meas_G$ is the complex vector $\Psi(\nu)\in\C^n$ with coordinates:
    \begin{equation*}
        \Psi(\nu)_l = \dualdot{\nu}{\psi_l}.
    \end{equation*}
\end{definition}

We consider three types of discrete transforms~$\Psi \colon \meas_G \to \C^n$. 

\paragraph*{Identity.}
The identity transform $\id \colon \meas_G \to \R^{ n}$ is associated to the family~$\{\psi_{i,j}\} = \{\1_{[b_i,b_{i+1})\times[p_j,p_{j+1})}\}$ induced by the finite sets $\mathcal{B}$ and $\mathcal{P}$. In that case, denoting the finite measure $\nu = \sum_{(i,j)\in G} m_{i,j}\delta_{(i,j)}$ for real numbers~$m_{i,j}$, we have $\Psi(\nu)_{i,j} = \dualdot{\nu}{\psi_{i,j}} = m_{i,j}$. Therefore, as its name suggests, the identity transform~$\Psi(\nu)$ is just another way of storing the finite measure $\nu$. We denote it by \qupid-\Id{} in our experiments.

\paragraph*{Fourier.}
The \emph{discrete Fourier transform} $\Fourier\colon \meas_G \to \C^{n}$ is associated to the family $\{\psi_{l_1,l_2}\}$ where $\psi_{l_1,l_2}(i,j) = \exp\left(-\mathbf{i}2\pi\left(l_1\frac{i}{r}+l_2\frac{j}{s}\right)\right)$; see \cite[Sec.~3.4]{mallat09}. This transform is computed using the Fast Fourier Transform algorithm~\cite{cooley1965algorithm}. The magnitude and phase of the output are flattened and concatenated to form the vectorization we use in our experiments, denoted by \qupid-\fft.


\paragraph*{Wavelets.}
The last type of transforms used in this article are discrete wavelet transforms, namely, the Daubechies and Coiflet transform; see~\cite[Chap.~7]{mallat09} and \cite{daubechies92}. The definition of wavelets on the two-dimensional plane is made precise in \cite[Thm.~7.25]{mallat09}. The \emph{Daubechies} and \emph{Coiflet wavelet transforms of order~$p$} are associated to families of wavelets satisfying specific theoretical properties. Namely, the Daubechies wavelets of order~$p$ have~$p$ vanishing moments~\cite[(7.69)]{mallat09} and have minimal support with this property~\cite[Thm.~7.7]{mallat09}. The Coiflet wavelets of order $p$ have minimal support with~$p$ vanishing moments and such that their so-called \emph{scaling function} \cite[Sec.~7.1.2]{mallat09} are normalized and have~$p-1$ vanishing moments.

The choice of the order $p$ yields different wavelets, hence allowing for further flexibility in our framework.
When~$p=1$, the Daubechies wavelet transform coincides with the simpler \emph{Haar wavelet transform}; see \cite[Sec.~7.2.2]{mallat09}. We provide the explicit definition of the Haar transform below. Defining the following two functions $\psi = \1_{[0,1/2)} - \1_{[1/2,1)}$ and~$\phi = \1_{[0,1)}$,
the Haar wavelet transform is associated to the family of wavelets $\{\psi^V_{k,l_1,l_2},\psi^H_{k,l_1,l_2},\psi^D_{k,l_1,l_2}\}_{(k,l_1,l_2)\in\Z^3}$ where:
\begin{align*}
    \psi^V_{k,l_1,l_2}(x_1,x_2) &= \frac{1}{2^k} \phi\left(\frac{x_1-2^kl_1}{2^k}\right)\psi\left(\frac{x_2-2^kl_2}{2^k}\right), \\
    \psi^H_{k,l_1,l_2}(x_1,x_2) &= \frac{1}{2^k} \psi\left(\frac{x_1-2^kl_1}{2^k}\right)\phi\left(\frac{x_2-2^kl_2}{2^k}\right), \\
    \psi^D_{k,l_1,l_2}(x_1,x_2) &= \frac{1}{2^k} \psi\left(\frac{x_1-2^kl_1}{2^k}\right)\psi\left(\frac{x_2-2^kl_2}{2^k}\right).
\end{align*}

In our experiments, we use the \texttt{Python} library \texttt{PyWavelets} to compute discrete transforms. For the Haar wavelets, taking as input a finite measure~$\nu\in\meas_G$ (such as  quantized persistence diagram) seen as a two-dimensional array, this library
outputs the vector $\Psi(\nu)$ of \Cref{def:transforms} associated to the family
\begin{equation}\label{eq:k-1-subfamily}
    \{\psi^V_{1,l_1,l_2},\psi^H_{1,l_1,l_2},\psi^D_{1,l_1,l_2}\}_{(l_1,l_2)\in\Z^2}
\end{equation}
obtained for $k=1$, called respectively \emph{vertical}, \emph{horizontal} and \emph{diagonal details}. In addition, it outputs the vector of \emph{approximation coefficients} obtained for the family:
\begin{equation*}
    \psi^A_{1,l_1,l_2}(x_1,x_2) = \frac{1}{2} \phi\left(\frac{x_1-2l_1}{2}\right)\phi\left(\frac{x_2-2l_2}{2}\right).
\end{equation*}
In practice, we concatenate horizontal, vertical, and diagonal details and approximation coefficients into a vectorization of the input measure $\nu$. We proceed similarly for the other Daubechies and Coiflet wavelets of order $p$, respectively denoted by \qupid-\db{p} and \qupid-\coif{p} in our experiments. 
Of course, one could get finer vectorizations by considering larger subfamilies than the one from~\eqref{eq:k-1-subfamily}, including other values of $k\in\Z$.

\subsection{Algorithmic complexity}\label{sec:complexity}
One of the main strengths of \qupid{} is its computational simplicity. Given a persistence diagram with~$l$ points, the quantization step consists solely of a binning process, hence will run in $\mathcal{O}(l)$ operations independently of the size $n = r \times s$ of the grid used for the quantization. The total time complexity of \qupid{} will therefore be of $\mathcal{O} (l+n \log (n))$ for the Fast Fourier Transform~\cite{cooley1965algorithm} and of $\mathcal{O}(l+n)$ for a Fast Wavelet Transform~\cite{beylkin1991fast}. It is to be noted that according to \cite{yogeshwaran2017random}, in some random settings, the number of points $l$ in the persistence diagram of a \v{C}ech filtration of a randomly sampled point cloud in the Euclidean space will asymptotically be of the same order as the number of points in the initial point cloud. 

This is to be paralleled by the computational complexity of computing persistence diagrams. For computational reasons, we do not use the \v{C}ech filtration in our experiments on point clouds but rather a homotopy equivalent one called the \emph{alpha filtration}; see \cite{bauer17}. Doing so, computing any degree persistence diagram from a point cloud with $N$ points in $\R^d$ has a time complexity of $\mathcal{O}(N^{d\omega/2})$ in the worst case, where $2 \leq \omega \leq 2.373$ is the matrix multiplication time; see~\cite{seidel95,milosavljevic2011zigzag}.

Taken together, this assesses that the computational cost of our vectorization is expected to be much smaller than the computation of persistence diagrams themselves. We will back up these theoretical claims by timing experiments in Section \ref{sec:orbit} by comparing our method with standard persistence diagrams vectorizations and showing that \qupid{} benefits from a reduced computational cost in practice.

\section{Experiments}
\label{sec:experiments}

In this section, we evaluate the vectorization procedure of \qupid{} by performing a classification task on some synthetic and real-world datasets. The typical workflow is the following: starting from a raw data input, we compute persistence diagrams using the \texttt{gudhi} library~\cite{maria2014gudhi}. We then apply \qupid{} for all possible transforms and a user-given grid resolution. Unless explicitly mentioned otherwise, the classification task itself is performed with a random forest classifier \cite{breiman2001random} using the implementation from the \texttt{scikit-learn} library~\cite{pedregosa2011scikit}. The code to reproduce the experiments is available at \href{https://github.com/majkevh/qupid}{\texttt{https://github.com/majkevh/qupid}}.

In what follows, we specify the transform used in every accuracy table or write \qupid{}$\star$ for the accuracy corresponding to the transform reaching maximum accuracy among all the ones tried.




\subsection{Dynamical system (\orbit5k dataset)}\label{sec:orbit}
\begin{figure}
    \centering
\includegraphics[width = \linewidth]{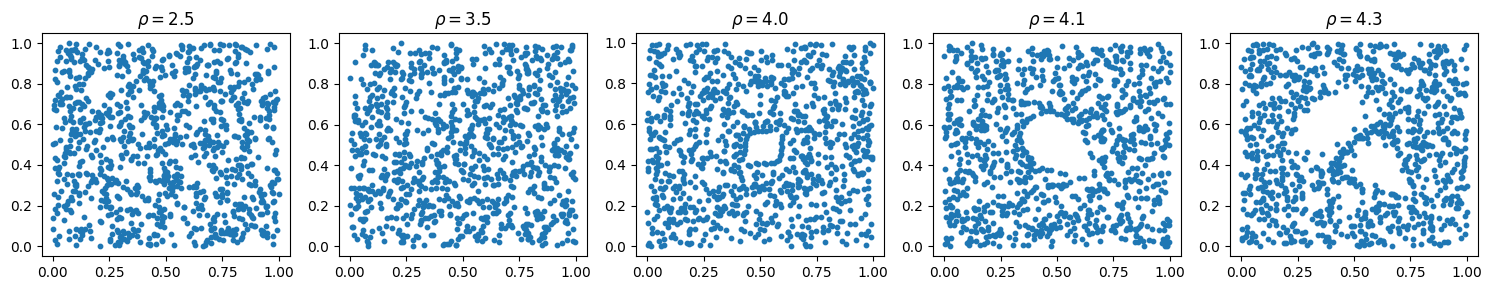}
    \caption{Examples of point clouds from the \texttt{ORBIT5K} dataset.}
    \label{fig:orbits}
\end{figure}

The \orbit5k dataset is a commonly used benchmark for supervised classification tasks in topological data analysis \cite{adams2017persistence, le2018persistence, carriere2017sliced, bba, pmlr-v108-carriere20a, pmlr-v130-royer21a,hacquard23}. It is a synthetic dataset consisting of point clouds of a thousand points in the unit square $[0,1]^2$, generated by a dynamical system depending on a real parameter~$\rho>0$. For each point cloud, the first point~$(x_0,y_0)$ is drawn uniformly at random in the unit square, and the remaining points $(x_n,y_n)$ for $n=1,\ldots,999$ are generated recursively by the formula:
\begin{align*}
    x_{n} & = x_{n-1} + \rho y_{n-1}(1 - y_{n-1}) &\mod 1, \\
    y_{n} & = y_{n-1} + \rho x_n(1 - x_n) &\mod 1.
\end{align*}
The whole dataset consists of five classes, each of a thousand point clouds, for five values of the parameter $\rho\in\{2.5,3.5,4.0,4.1,4.3\}$. See \Cref{fig:orbits} for an illustration.

\paragraph*{Supervised classification.}

The dataset is divided by a 70/30 randomized train/test split for each class. Changing the parameter $\rho$ changes the underlying topology and sampling distribution of the generated orbits. We try to predict the value of $\rho$ given the point cloud as data. We compare two vectorizations based on our methodology described in \Cref{sec:methodology}. We begin by computing degree-0 and degree-1 persistence diagrams of the \v{C}ech filtration computed on the point clouds as described in \Cref{sec:background}. As mentioned in \Cref{sec:complexity}, we use alpha-filtrations instead of \v{C}ech filtrations for computational reasons; see \cite{bauer17}. We consider an adaptive log-scaled grid with a parameter $\alpha$ chosen by grid search cross-validation, coupled with the Fourier transform and the Daubechies and Coiflet wavelet transforms of orders $1$, $2$ and $3$. We also considered the identity transform over a uniform grid, noticing that the choice of an adaptive grid gave highly sub-optimal results. Our accuracies are averaged over ten runs and reported for the transform reaching higher accuracy (Coiflet of order 2 in this case) in \Cref{tab:orbit}. Additional results for other discrete transforms are available in the appendix (\Cref{tab:complete:orbit}).

In \Cref{tab:orbit}, we also compare our method with three state-of-the-art vectorizations of persistence diagrams. Namely, Persistence Images (PI) \cite{adams2017persistence}, Persistence Fisher Kernels (PF-K) \cite{le2018persistence} and Persformer \cite{reinauer2021persformer}.

\begin{table}[ht]
    \centering
    \resizebox{\linewidth}{!}{
    \begin{tabular}{|c|c|c||c|c|}
    \hline PI$^\star$ & PF-K$^\dagger$ & Persformer$^\ast$ & \qupid$\star^\star$ & \qupid{}-\Id$^\star$ \\
    \hline $82.5$ & $85.9\pm.8$  & $91.2\pm.8$ & $88.3\pm.4$ & $82.0\pm.6$\\
    \hline
    \end{tabular}}
    \caption{Classification scores for the \orbit5k dataset. Notes: $^\dagger$ denotes methods with accuracy averaged over 100 runs; $^\star$ indicates averaging over 10 runs; $^\ast$ signifies an average over 5 runs. No standard deviation is reported in \cite{adams2017persistence} for persistence images.}\label{tab:orbit}
\end{table}

On this classification task, taking wavelet transforms induces a significant gain in accuracy allowing \qupid{} to reach a state-of-the-art accuracy. It is only outperformed by the Persformer method which relies on the costly training of a neural network with transformer architecture. In comparison, our method is faster to compute than all other methods, as shown in the following paragraph. 

Finally, let us also point out that applying \qupid-\Id{} directly on the point cloud without any homology computation yields an accuracy of around~$99 \%$, similar to what was observed in \cite{pmlr-v130-royer21a, reinauer2021persformer}. Therefore, this experiment should be considered only as a proof-of-concept to benchmark our method against other vectorizations of persistence diagrams.

\paragraph*{Timing.} All computations were conducted on a single laptop, an Apple M3 Pro equipped with an 11-core CPU and no GPU. The main strength of our method lies in its simplicity, which translates into a very meagre computational cost while still keeping a competitive accuracy. This is shown in \Cref{tab:times} where we report the time to vectorize the \orbit5k dataset for all methods cited in \Cref{tab:orbit}. We use a resolution of~$32\times 32$ for persistence images (PI) and \qupid{} with log-scaled grids, and the PF-K vectorization is timed with default parameters. Computational times for Persformer are not provided in~\cite{reinauer2021persformer} and training this Transformer architecture on our CPU was impossible in a reasonable time. Still, the authors of the paper have shared with us in a private communication that the training process takes in the order of hours on GPUs, which we report in the table.
\begin{table}[h]
    \centering
    \begin{tabular}{lc}
        \toprule
        {Vectorization} & {Time} \\
        \midrule
        \texttt{\qupid-\Id} & 15 s \\
        \texttt{\qupid-\fft} & 17 s \\
        \texttt{\qupid-\db{2}} & 15 s \\
        \texttt{\qupid-\coif{2}} & 16 s \\
        PF-K & 69 s \\
        PI & 74 s  \\
        Persformer & $\sim$ hours \\
        \bottomrule
    \end{tabular}
    \caption{Time to vectorize the \orbit5k dataset.}
    \label{tab:times}
\end{table}

First, note that considering discrete Fourier or wavelet transforms has a negligible cost compared to the quantization process itself. Overall, the results of Tables \ref{tab:orbit} and \ref{tab:times} show the excellent accuracy-complexity trade-off of \qupid. On the one hand, our method allows for an increase of almost 3\% in accuracy with a computational time divided by more than 4 for PI and PF-K. On the other hand, our method takes only a dozen of seconds on a CPU against hours for Persformer on a GPU at the cost of only 3\% in accuracy.

\subsection{Graph classification}\label{sec:graphs}
All the previous examples dealt with point clouds as input data. Topological data analysis can take a much larger class of input data, such as graphs. In this section, we consider a supervised graph classification problem. To compute the persistence diagrams of a graph, we need to introduce a topological filtration. A common choice initiated in \cite{pmlr-v108-carriere20a} is to introduce the heat kernel signature of a graph:
\begin{definition}\label{def:hks}
  For a graph $\mathcal{G}=(V, E)$, the Heat Kernel Signature (HKS) function with diffusion parameter $t > 0$ is defined for each vertex $v \in V$ as
$$
\HKS_t(v)=\sum_{k=1}^{|V|} \exp \left(-t \lambda_k\right) \psi_k(v)^2,
$$
where $\lambda_k$ is the $k$-th eigenvalue of the normalized graph Laplacian, and $\psi_k$ is the corresponding eigenfunction.
\end{definition}

The heat kernel signature at a vertex $v$ can be interpreted as the remaining heat at $v$ starting from a unit heat distribution on $v$ after a time $t$ has elapsed; see~\cite{sun2009concise}. The HKS carries a lot of local and global structural information on the graph, related to its connectivity and topology. The HKS defines a function on the set of vertices that can be simply extended to the set of edges. If $e=(x,y)$ is an edge between vertices $x$ and $y$, we set $\HKS_t (e) = \max(\HKS_t(x), \HKS_t(y))$. For a given parameter~$t$, considering the graphs given by the sublevel sets or superlevel sets of the HKS gives us a topological filtration. 

We adopt the same methodology as in \cite{pmlr-v108-carriere20a} and \cite{pmlr-v130-royer21a}, where taking the degree-0 and degree-1 persistence diagrams both for the sublevel sets and superlevel sets of the HKS provides us with a family of four persistence diagrams per graph called \textit{extended persistence}.

We benchmark our method with state-of-the-art graph classification algorithms on a few real-world datasets. The datasets include the social graphs \texttt{IMDB-B}, \texttt{IMDB-M}, \texttt{COLLAB}, \texttt{REDDIT-5K} and \texttt{REDDIT-12K}, as well as graphs coming from chemoinformatics and bioinformatics, namely \texttt{COX2}, \texttt{DHFR}, \texttt{MUTAG} and \texttt{PROTEINS}. 

 We adopt the methodology from  \cite{pmlr-v130-royer21a} for generating extended persistence diagrams, using two specific HKS diffusion times ($t_1=0.1$ and $t_2=10$), focusing solely on the topological structure without considering additional graph attributes. Therefore, we end up with eight persistence diagrams and concatenate all corresponding vectorizations obtained from \qupid{}.
For the evaluation, we used three distinct regular grid sizes~$20\times20$, $32\times32$, and $50\times50$ across all our proposed transforms except for the \texttt{COLLAB} and \texttt{REDDIT}s datasets where we used only~$10\times10$ and $20\times20$. 
Regarding the wavelet transforms, we experimented with six different wavelet types, namely Daubechies and Coiflet of orders~$1$, $2$ and~$3$. 
The highest accuracy achieved over these transforms (\qupid{}$\star$) is presented in Table \ref{results:graph}. Detailed accuracy for each transform across the different grids and standard deviations can be found in the appendix (\Cref{tab:complete:graphs}). Performance is assessed through 10-fold evaluations on each dataset, reporting both average and best 10-fold results.
\begin{table*}[t]
    \centering
    \resizebox{1\textwidth}{!}{
    \begin{tabular}{|l||ccccc||ccc||cc|}
    \hline 
     & SV$^\dagger$ & RetGK$^{\star}$ & FGSD$^\dagger$ & GCNN$^\dagger$ & GIN$^\dagger$ & P$^\dagger$ & Perslay$^\star$ &ATOL$^\star$ & \multicolumn{2}{c|}{\qupid{}$\star$} \\
    Dataset & &&&& & && &Mean$^\star$ & Max$^\dagger$ \\
    \hline 
    \texttt{MUTAG} & $88.3$ & $90.3$ & $92.1$ & $86.7$ & $89.0$ & $79.2$ & $89.8 $ & $88.3  $ & $89.8$ & $91.7$ \\
    \texttt{COX2} & $78.4$ & $81.4$ & - & - & - &$76.0$ & $80.9$ & $79.4$ & $80.6$ & $82.4 $\\
    \texttt{PROTEINS} & $72.6$ & $75.8$ & $73.4$  &  $76.3$ & $75.8$ &$65.4$ &$74.8$ & $71.4$ &$72.8$ &$ 73.6 $\\
    \texttt{DHFR}  & $78.4$ & $81.5$ & - & - & - & $70.9$ & $80.3$ &$82.7$ & $81.8$ & $83.1$\\
    \texttt{IMDB-B} & $72.9$ & $71.9$ & $73.6$ & $73.1$ & $74.3$ & $54.0$ & $71.2$ & $74.8$ & $68.2$& $69.3$ \\
    \texttt{IMDB-M} & $50.3$ & $47.7$ & $52.4$ & $50.3$ & $52.1$ & $36.3$ & $48.8$& $47.8$ & $43.7$& $44.4$\\
    \texttt{COLLAB} & - & 81.0 & 80.0 & 79.6 & 80.1 & - &76.4 & 88.3 & 86.4& 86.7   \\
    \texttt{REDDIT-5K} & - & 56.1 & 47.8 & 52.9 &   57.0 & - & 55.6 & 67.1 & 66.1 & 66.6 \\
    \texttt{REDDIT-12K} & - & 48.7 & - &  46.6& - &  - & 47.7 & 51.4 & 51.4& 51.7\\
    \hline
    \end{tabular}}
    \caption{Classification scores for graph datasets. Note: The symbol $^\dagger$ indicates that the performance should be compared with the \textit{max} metric, while the symbol $\star$ denotes that the comparison should be made with the \textit{mean} metric.}\label{results:graph}
\end{table*}

In \Cref{results:graph}, we also compare ourselves with state-of-the-art graph classification methods. In the first three columns, we compare with methods using neural networks or kernel methods. Namely, a scale-variant topological method (SV from \cite{tran2018scalevariant}), a graph kernel method based on random walks (RetGK from \cite{zhang2018retgk}), a graph embedding method using spectral distances (FGSD from \cite{verma2017hunt}),  a Capsule Graph Neural Network (GCNN from \cite{zhang2019capsule}), and a Graph Isomorphism Network (GIN from \cite{xu2019powerful}). In the three following columns, we compare our results with Persistence Images from \cite{adams2017persistence} and Persistence Landscapes from \cite{bubenik2015statistical} coupled with an XGboost classifier. The maximum accuracy of these two methods is denoted in the ``P'' column. We also compare ourselves with Perslay~\cite{pmlr-v108-carriere20a} and ATOL~\cite{pmlr-v130-royer21a}.

These results show that despite its apparent simplicity, \qupid{} achieves decent performances in all classification tasks. We always outperform persistence images and landscapes (P) and are comparable with other TDA-based vectorizations of persistence diagrams across all datasets. Importantly, we show very competitive performances on the \texttt{COLLAB} and \texttt{REDDIT}s datasets, comparable to the ATOL algorithm, thus outperforming all other methods. This should be read in the light of the fact that, as explained in~\cite{pmlr-v130-royer21a}, those three datasets are considered the most difficult to classify due to their size in terms of number of data and number of edges and vertices.

\subsection{Tumor immune cells}\label{sec:cancer}
We use our method to study a real-world dataset introduced in \cite{vipond21} which contains point clouds in $\mathbb{R}^2$ obtained from spatial locations of immune cells in tumours. The dataset contains three different types of immune cells, namely, \CDs (65 point clouds with average of 560 points, minimum of 88 and maximum of 1000), \Fox (74 point clouds with average of 211 points, minimum of 24 and maximum of 499) and \CDh (73 point clouds with average of 651 points, minimum of 194 and maximum of 1000). In~\cite{vipond21}, authors use the sophisticated machinery of \emph{multiparameter persistent landscapes} to study this dataset for a multi-parameter filtration, namely, position and local density of cells. Here, we study the same data using only one parameter, following our methodology introduced in \Cref{sec:methodology}, and used in \Cref{sec:orbit}. We compute degree-0 and degree-1 persistence diagrams associated with the \v{C}ech filtration (in fact, alpha filtration) and consider the \qupid{} vectorization.

\begin{figure}[H]
    \centering
\includegraphics[width = 1\linewidth]{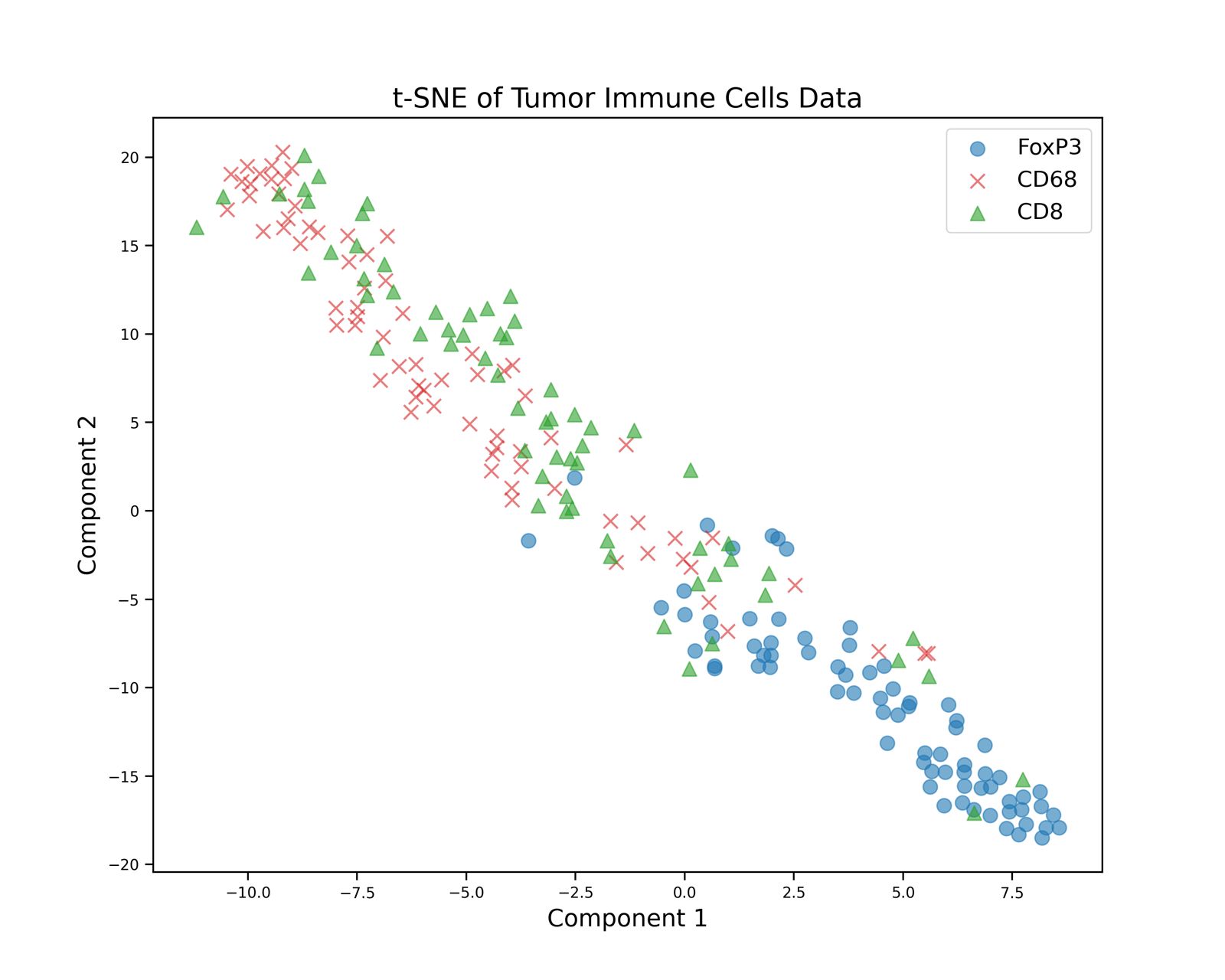}
    \caption{t-SNE plot of \qupid{} for the Coiflet transform of order $2$ on the tumor immune cells dataset.}
    \label{fig:pca:cancer}
\end{figure}

\paragraph*{Unsupervised setting.} 
We start by considering the problem of unsupervisedly predicting the type of immune cell. We show in \Cref{fig:pca:cancer} a t-SNE plot of \qupid{} using the Coiflet transform of order $2$ on the degree-0 and degree-1 persistence diagrams. Our method satisfyingly separates the \Fox class from the two other classes suggesting that this class is easy to discriminate from the other two. This shows that \qupid{} already highlights discriminative information on this dataset at the unsupervised level. 

\paragraph*{Supervised setting.}
We use our methodology in the supervised cell-type classification task and the three corresponding binary problems.  We use the same setting as in \Cref{sec:orbit} to vectorize persistence diagrams using \qupid{}, with adaptive log-scaled grids. We denote it by \qupid{}$_\mathrm{PD}$. For each classification problem, we make a randomized 80/20 training/test split as in \cite{vipond21}. We report the results averaged over 20 runs in \Cref{tab:cancer}. Detailed accuracies for each transform across the different transforms can be found in the appendix~(\Cref{tab:complete:cancer}).

We also compare our vectorization of persistence diagrams with four other methods. First, we compare ourselves to scores reported in the original paper~\cite{vipond21} using multiparameter persistence landscapes (MPL), a topological descriptor using a two-parameter filtration. Second, we vectorize the raw point clouds using \qupid{} with a uniform grid of size $32\times 32$, comparing the same eight transforms: identity, Fourier, Daubechies ($p=1,2,3$) and Coiflet ($p=1,2,3$). We denote it by~\qupid$_\mathrm{raw}$. The last two methods consist in using ATOL~\cite{pmlr-v130-royer21a} to vectorize the raw point clouds (ATOL$_\mathrm{raw}$) or their persistence diagrams (ATOL$_\mathrm{PD}$). We chose the ATOL budget-parameter giving the best accuracy. For every vectorization, we applied a random forest or an LDA classifier from the \texttt{scikit-learn} library and reported the best accuracy for each method in \Cref{tab:cancer}.
\begin{table}[H]
\centering
\adjustbox{max width=\linewidth}{
\begin{tabular}{|c|c|c|c||c|c|c|}
\hline
Task & MPL & ATOL$_\mathrm{PD}$ & ATOL$_\mathrm{raw}$&\qupid$_\mathrm{PD}\star$ & \qupid$_\mathrm{raw}\star$\\
\hline
$\texttt{CD68}$ vs $\texttt{CD8}$ &65.3& 61.1&51.4& \textbf{79.0} &65.0\\
$\texttt{CD68}$ vs $\texttt{FoxP3}$ &86.3&90.7 & 87.7&\textbf{98.2} & 89.2\\
$\texttt{CD8}$ vs $\texttt{FoxP3}$ &74.7&85.9 &85.7&\textbf{96.3} & 87.7\\
All classes & 66.2 & 67.3 &58.5& \textbf{79.4}&67.9\\
\hline
\end{tabular}}
\caption{Comparative analysis of performance for immune cells classification task.}
\label{tab:cancer}
\end{table}

First, results of \Cref{tab:cancer} demonstrate the advantage of considering persistence diagrams over raw point clouds, both for the \qupid{} and the ATOL vectorizations. This suggests that the topological and geometric information captured by persistence diagrams is relevant for this classification task.

Quantitatively, these results demonstrate the benefits of using the \qupid{} algorithm to vectorize persistence diagrams. Not only do we outperform ATOL vectorization of persistence diagrams in all tasks---and by more than 10\% in all but one---, but we also outperform the multiparameter approach MPL by even larger amounts---22\% in the case of \CDh{} vs \Fox{}. We believe this is due to the ability of our method to focus on low-persistence points, allowing it to capture fine geometric information missed by other methods.

\paragraph*{Timing.} Finally, we report the computational times for all methods in this classification task. Computing multiparameter persistence landscapes on the entire dataset took approximately 3 hours, against a few seconds for \qupid{} and ATOL. Moreover, computing persistence diagrams is the most time-consuming task, with more than 1 second, against less than 0.5 seconds for the vectorization step.

\begin{figure*}[t]
    \centering
    \includegraphics[width=\linewidth]{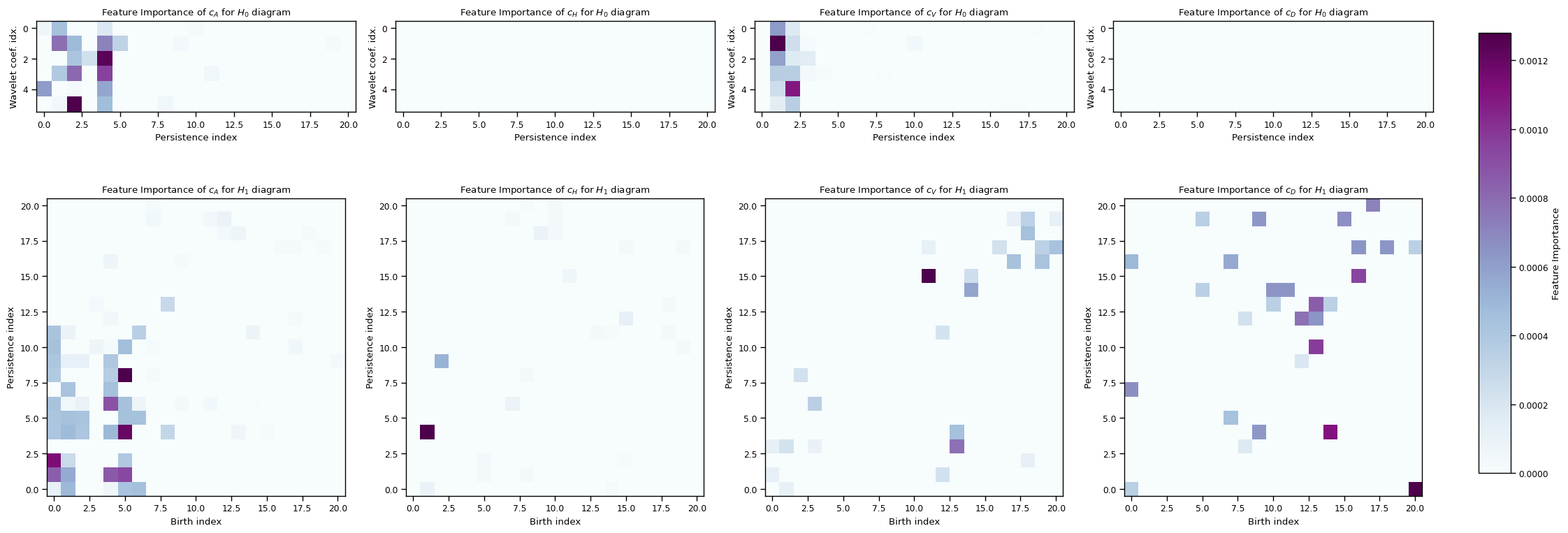}
    \caption{Importance of the coefficients of the \qupid-\coif{2} vectorization for the classification task \CDs{} vs \Fox{}.}
    \label{fig:feature-importance}
\end{figure*}

\paragraph*{Feature importance.} To gain some interpretability, we compute the random forest's feature importance of each coefficient of the vectorizations using \texttt{scikit-learn}. In \Cref{fig:feature-importance}, we plot the feature importance heatmap of the \qupid-\coif{2} vectorization of persistence diagrams for the binary classification of \CDs{} vs \Fox{}. We can see that both the degree-0 and degree-1 persistence diagrams are used in the classification, but the former is only through its approximation coefficients and vertical detail. For degree-1 persistence, the approximation coefficients carry the most information, mainly through its early-birth features. The horizontal and vertical details also contain important information, extremely concentrated in a few coefficients. The possibility to make such an analysis is crucial in sensitive applications such as processing medical data, where we can directly relate the prediction to a few key-details in the wavelet analysis of topological descriptors.

\subsection{Influence of the resolution parameter}

\begin{figure}
    \centering
    \includegraphics[width = \linewidth]{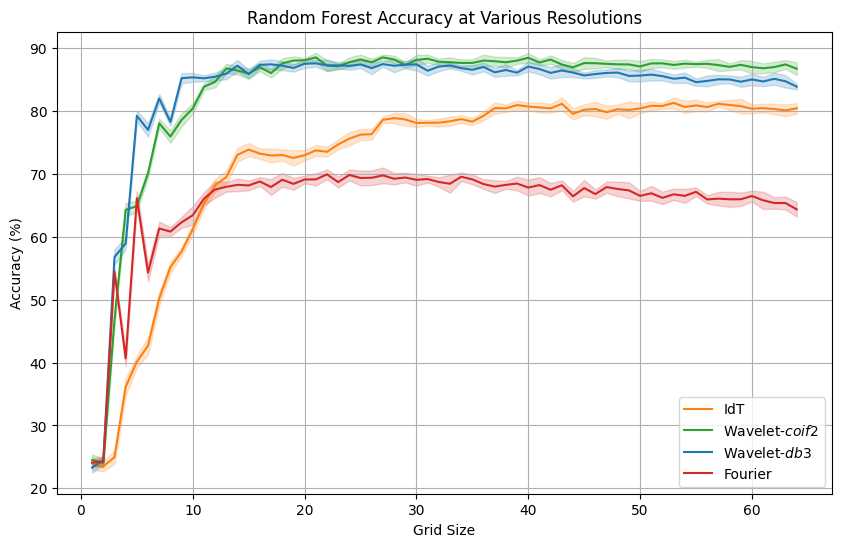}
    \caption{Influence of the resolution parameter on the accuracy for classifying the \orbit5k dataset.}
    \label{fig:ablation}
\end{figure}

As demonstrated in the previous quantitative experiments, considering discrete transforms usually yields higher accuracy than using the identity transform in supervised tasks. To further stress this fact, we come back to the setting of Section \ref{sec:orbit} where we compare the accuracy of several transforms on a square grid of varying size; see \Cref{fig:ablation}. A similar study has been conducted in~\cite{pmlr-v130-royer21a,hacquard23}.

Naturally, increasing the resolution usually provides a higher accuracy. We notice that for this task, wavelet transforms offer better accuracy than the identity at every possible resolution and that the difference is even more striking at low resolutions. Indeed, benefiting from their signal-compression properties, wavelet transforms can offer a competitive accuracy using a very limited number of coefficients. This makes them a first-choice option when working under a limiting vectorization budget. 

\section{Conclusion and further work}\label{sec:conclusions}
This work falls within an old endeavour of mapping persistence diagrams into a vector space to use off-the-shelf statistical learning methods. We believe that the algorithm developed in this article allows us to take a step back from the existing literature and to realize that using a seemingly naive approach of binning the diagrams and using standard discrete transforms accounts for strong predictive power. More precisely, the take-home message of the \qupid{} algorithm is the following:
\begin{itemize}
    \item We reach a very good accuracy in a wide variety of settings.
    \item The gap in accuracy with the most competitive methods is traded off by the extreme simplicity of the method and its reduced computational time.
    \item We successfully made use of log-scaled grids that focus on the low-persistence points, giving one more piece of evidence that the near-diagonal points carry a lot of discriminative information.
    \item In addition to improving the overall accuracy, using Fourier and wavelet transforms allows for a certain robustness to a resolution loss, successfully compressing the topological information. 
    \item By its simplicity, we claim that we have developed a white-box algorithm allowing for immediate interpretability. This property becomes extremely valuable when working on real-world data.
    \end{itemize}

Using discrete transforms to analyze topological descriptors is a research direction of growing interest, both theoretically; see \cite{haberle2024wavelet} and experimentally; see \cite{hacquard23}. In the present paper, we have benefitted from the meagre computational cost of our algorithm to optimize over the best integral transform and grid parameters. Future works involving finer qualitative approaches to better understand which transform is more suited to a given type of data would be essential. Finally, we point out that \qupid{} is perfectly suited for the vectorization of many recent invariants of multiparameter persistent homology, which take the form of discrete measures in the plane~\cite{loiseaux23}. We defer the study of \qupid{} in such multiparameter analyses for future work.

\bibliographystyle{unsrt}
\bibliography{reference} 

\newpage 
\appendix

\onecolumn
\section{Complementary material}\label{appendix:B}
This section includes extensive tables that provide a comparative analysis between every transform on the classification tasks of Section \ref{sec:experiments}. 

\begin{table*}[h]
\centering
\adjustbox{max width=\linewidth}{
\begin{tabular}{|c|c|c|c||c|c||c|c||c|c||c|c||c|c||c|c||c|c|}
\hline
\multicolumn{2}{|c|}{Task \multirow{2}{*}{}} & \multicolumn{2}{c||}{\qupid-\fft{}} & \multicolumn{2}{c||}{\qupid-\db{1}} & \multicolumn{2}{c||}{\qupid-\db{2}} & \multicolumn{2}{c||}{\qupid-\db{3}} & \multicolumn{2}{c||}{\qupid-\coif{1}} & \multicolumn{2}{c||}{\qupid-\coif{2}} & \multicolumn{2}{c||}{\qupid-\coif{3}} & \multicolumn{2}{c|}{\qupid-\Id{}} \\
\cline{1-18}
Dataset & Grid size & Mean & Max & Mean & Max & Mean & Max & Mean & Max & Mean & Max & Mean & Max & Mean & Max & Mean & Max \\
\hline
\multirow{3}{*}{\texttt{MUTAG}} & $20\times 20$ & $87.8(\pm.8)$ & $88.9$ & $87.8(\pm.8)$ & $88.9$ & $87.8(\pm.8)$ & $88.9$ & $88.7(\pm1.0)$ & $90.6$ & $88.6(\pm.8)$ & $90.0$ & $87.7(\pm1.2)$ & $90.0$ & $87.5(\pm.9)$ & $89.4$ & $88.1(\pm.7)$ & $89.4$ \\
& $32\times 32$ & $89.1(\pm1.2)$ & $91.1$ & $88.4(\pm1.3)$ & $90.0$ & $88.6(\pm1.0)$ & $90.0$ & $88.9(\pm1.0)$ & $90.6$ & $88.9(\pm.7)$ & $90.0$ & $89.7(\pm.5)$ & $90.6$ & $\mathbf{89.8(\pm.9)}$ & $\mathbf{91.7}$ & $89.1(\pm.9)$ & $90.7$ \\
& $50\times 50$ & $88.1(\pm1.1)$ & $89.4$ & $86.9(\pm1.1)$ & $88.9$ & $87.6(\pm.9)$ & $88.9$ & $88.3(\pm1.4)$ & $91.1$ & $88.2(\pm.9)$ & $89.4$ & $88.3(\pm1.1)$ & $90.0$ & $88.1 (\pm.6)$  & $88.9$ & $87.4(\pm1)$ & $88.9$ \\
\hline
\multirow{3}{*}{\texttt{DHFR}} & $20\times 20$ & $77.8(\pm.7)$ &$78.9$ & $77.8(\pm.6)$ &$78.9$ & $78.3(\pm.7)$& $78.9$ &$79.3(\pm.6)$ & $80.0$&$78.5(\pm.6)$ & $79.6$&$79.2(\pm1.0)$ &$81.2$ & $79.0(\pm.8)$& $80.4$&$78.7(\pm.9)$ &$79.7$ \\
& $32\times 32$ &$79.1(\pm1.0)$ &$80.5$ & $81.5(\pm.7)$&$\mathbf{83.1}$ &$\mathbf{81.8(\pm.7)}$ & $82.9$& $81.1(\pm1.0)$ & $82.4$& $81.4(\pm.5)$ & $82.1$& $80.7(\pm.5)$&$82.4$ & $80.6(\pm.7)$& $81.9$ & $81.1(\pm.5)$& $82.0$ \\
& $50\times 50$ &$79.6(\pm.7)$ &$81.2$ & $ 80.6(\pm.4)$& $81.2$& $80.2(\pm.4)$&$80.8$ &$80.7(\pm.4)$ &$81.5$ & $81.0(\pm.5)$&$81.9$ & $81.5(\pm.5)$&$82.3$ & $80.9(\pm.9)$& $82.7$&$81.2(\pm.7)$& $82.7$ \\
\hline
\multirow{3}{*}{\texttt{COX2}} & $20\times 20$ &$78.3(\pm.6)$ &$78.9$ & $79.2(\pm.7)$& $80.4$& $79.7(\pm.7)$& $81.1$& $79.7(\pm.7)$&$80.9$ & $80.0(\pm.4)$&$80.4$ & $79.7(\pm.5)$&$80.4$ &$79.5(\pm.8)$ & $81.5$& $79.4(\pm.7)$& $80.4$ \\
& $32\times 32$ & $76.7(\pm.7)$ & $77.8$& $77.3(\pm1.0)$&$79.3$ & $78.0(\pm1.0)$&$79.3$ & $78.0(\pm1.1)$&$79.3$ & $78.0(\pm.6)$& $78.9$& $78.2(\pm.6)$& $78.9$& $77.4(\pm.8)$& $78.5$& $77.5(\pm.5)$& $78.5$ \\
& $50\times 50$ &$77.1(\pm.9)$ & $78.3$ & $77.5(\pm.7)$ & $78.5$& $78.2(\pm.5)$&$78.9$ &$77.7(\pm.6)$ & $78.9$& $78.7(\pm.7)$&$80.0$ & $78.2(\pm1.0)$&$79.3$ &$78.8(\pm1.2)$ & $80.2$& $78.0(\pm.4)$& $78.5$ \\
\hline
\multirow{3}{*}{\texttt{PROTEINS}} & $20\times 20$ & $72.5 (\pm.6)$  & $73.2$ & $72.4 (\pm .4)$ & $72.8$&$72.5 (\pm .7)$  &$73.4$ & $72.6 (\pm .6)$ & $73.4$  & $\mathbf{72.8 (\pm .6)}$ & $\mathbf{73.6}$ & $72.3 (\pm .5) $& $73.0$ & $72.1 (\pm.4)$ & $72.9$& $71.9(\pm.7)$& $73.2$ \\
& $32\times 32$ & $71.9 (\pm .5) $ & $72.8$ & $72.1 (\pm .5)$ & $72.9$ & $72.3 (\pm .5)$ & $73.1$ & $72.3 (\pm .5)$ & $73.3$ & $72.3 (\pm .6)$ & $73.4$ & $71.9 (\pm.5)$& $73.0$& $72.1 (\pm .8) $&$73.0$& $72.2(\pm.5)$& $73.2$ \\
& $50\times 50$ & $71.9(\pm .6$ & $72.7$ &$72.1 (\pm .5)$ &$73.1$ &$72.4 (\pm .6)$ &$73.2$ & $72.4 (\pm .6)$&$73.2$ & $72.7 (\pm.5)$&$73.5$ &$72.6 (\pm .5)$ &$73.5$ &$72.6 (\pm .5)$ & $73.3$& $71.4(\pm.4)$& $72.1$ \\
\hline
\multirow{3}{*}{\texttt{IMDB-B}} & $20\times 20$ &$67.6 (\pm.6)$ &$68.6$  &$67.5 (\pm.8)$ &$68.7$ & $67.1 (\pm .6)$& $68.3$& $67.0(\pm.6)$&$67.7$ & $67.1 (\pm.3)$&$67.5$ &$66.8 (\pm.7)$ &$67.9$ &$66.4 (\pm.9)$ & $68.3$ &$67.2(\pm.7)$& $68.3$ \\
& $32\times 32$ & $66.7 (\pm .5)$&$67.7$  &$66.4 (\pm.9)$ &$68.3$ &$66.9 (\pm.8)$ & $67.9$&$67.2 (\pm .6)$ &$67.9$ &$66.7 (\pm .8)$ &$67.7$ &$66.8 (\pm .6)$ &$67.7$ &$66.8 (\pm .6)$ &$67.8$ & $66.5(\pm.6)$&  $67.9$\\
& $50\times 50$ & $67.5 (\pm.8)$&$69.1$  &$\mathbf{68.2 (\pm.6)}$ &$\mathbf{69.3}$ &$68.1(\pm .7)$ & $69.2$& $67.4 (\pm.8)$ &$68.5$ &$68.2(\pm .7)$ &$69.1$ &$67.7 (\pm .4)$ &$68.3$ &$67.7 (\pm.5)$ &$68.5$&$66.7(\pm.9)$& $67.9$ \\
\hline
\multirow{3}{*}{\texttt{IMDB-M}} & $20\times 20$ &$42.2 (\pm .5)$ &$43.1$  &$42.4 (\pm.4)$ &$42.9$ &$42.2(\pm.5)$ &$42.9$ &$42.5 (\pm.3)$ &$43.3$ &$42.3 (\pm .7)$ &$43.5$ &$42.3 (\pm .3)$ &$42.8$ &$42.5 (\pm.4)$ &$43.2$ & $42.8(\pm.4)$& $43.6$ \\
& $32\times 32$ &$42.8 (\pm.5)$ &$43.7$  & $42.8 (\pm.5)$&$43.9$ &$42.5 (\pm.4)$ &$43.5$ &$43.1 (\pm.5)$ &$43.9$ &$42.9(\pm .5)$ &$43.6$ &$42.7 (\pm .4)$ &$43.5$ & $42.6 (\pm.5)$& $43.8$ & $43.1(\pm.5)$& $44.3$ \\
& $50\times 50$ & $42.9 (\pm.5)$ &$43.7$ & $42.9 (\pm.2)$ &$43.3$ &$43.0 (\pm.5)$ & $43.6$& 
 $43.3 (\pm .6) $&$\mathbf{44.4}$ &$43.1 (\pm .5)$ & $43.9$  &$\mathbf{43.7(\pm .5)} $ & $43.9$ & $43.3 (\pm .4)$ & $44.1$ & $43.1(\pm.4)$& $43.7$ \\
\hline
 \multirow{2}{*}{\texttt{COLLAB}} & $10\times 10$ & $86.1 (\pm.9)$ &$86.3$ & $85.7 (\pm.9)$&$85.9$ & $86.0 (\pm.2)$& $86.3$& $86.0 (\pm.6)$& $86.3$&$86.0 (\pm.6)$ & 
 $86.2$&$86.1 (\pm.2)$ &$86.3$ & $85.9 (\pm.2)$ & $86.3$& $85.7 (\pm.1)$&$85.9$ \\
 & $20\times 20$&$\mathbf{86.4 (\pm1.1)}$ & $86.5$ &$86.0 (\pm.2)$& $86.3$&$86.2 (\pm.2)$ & $86.4$& $86.3 (\pm.2)$&$86.6$ & $86.3 (\pm.1)$& $86.5$& $86.3 (\pm.1)$& $\mathbf{86.7}$& $86.3 (\pm .2)$ & $86.6$& $86.0 (\pm.2)$&$86.3$ \\
\hline
 \multirow{2}{*}{\texttt{REDDIT-5K}} &$10\times 10$  & $65.6(\pm.2)$&$65.9$ & $65.8(\pm.1)$ & $65.9$&  $\mathbf{66.1(\pm.2)}$& $66.5$&$65.8(\pm.3)$ & $66.2$&$65.7(\pm.3)$ & $66.2$&$65.6(\pm.3)$ &$66.2$ & $65.6(\pm.2)$& $65.9$& $65.5 (\pm.3)$& $66.1$\\
 & $20\times 20$ & $65.8 (\pm.2)$&$66.2$ & $65.9(\pm.2)$ & $66.2$& $66.0 (\pm.3)$ &$\mathbf{66.6}$ & $65.9 (\pm.3)$& $66.4$&$65.6 (\pm.2)$ & $66.0$&$66.0(\pm.3)$ &$66.5$ &$65.6(\pm.2)$ & $65.8$& $65.2(\pm.2)$& $65.5$\\
\hline
 \multirow{2}{*}{\texttt{REDDIT-12K}} & $10\times 10$ & $\mathbf{51.4 (
 \pm .1)}$ & $\mathbf{51.7}$& $50.8 (\pm .1)$ &$51.1$ & $51.1 (\pm.2)$& $51.6$& $51.3 (\pm.2)$&$51.6$ &$51.2(\pm.2)$ &$51.8$ & $51.0(\pm.2)$&$51.4$ & $51.0 (\pm.2)$&$51.3$ & $49.1 (\pm.3)$ & $49.8$\\
 & $20\times 20$& $51.0 (\pm.2)$ & $51.3$& $49.0 (\pm.2)$ & $(49.6)$& $50.1 (\pm.3)$& $50.5$& $50.5(\pm.2)$& $50.8$& $50.4(\pm.2)$& $50.8$&$50.3(\pm.2)$ &$50.7$ &$50.0(\pm.3)$ & $50.4$& $47.0 (\pm.2)$&$47.4$ \\
\hline
\end{tabular}
}
\caption{Comparative analysis of graph data classification scores of the methods presented in this study.}
\label{tab:complete:graphs}
\end{table*}

\begin{table*}[h]
\centering
\adjustbox{max width=\linewidth}{
\begin{tabular}{|c|c||c||c||c||c||c||c||c|c|}
\hline
Task & \qupid-\fft{} &\qupid-\db{1} & \qupid-\db{2} & \qupid-\db{3} & \qupid-\coif{1} & \qupid-\coif{2} &\qupid-\coif{3} &\texttt{\qupid-\Id{}} \\
\hline

\orbit5k  & $69.2 (\pm .7)$ & $82.0 (\pm1.0)$ & $87.9 (\pm.8)$ & $88.2(\pm.7)$ &$87.9 (\pm .6)$ & $\mathbf{88.3(\pm.4)}$ & $87.6(\pm.3)$ & $82.0\pm.6$\\
\hline
\end{tabular}}
\caption{Comparative analysis of \orbit5k classification scores for several choice of discrete transforms.}\label{tab:complete:orbit}
\end{table*}


\begin{table*}[h]
\centering
\adjustbox{max width=\linewidth}{
\begin{tabular}{|c|c|c|c|c|c|c|c|c|c|c|c|c|c|c|c|c|c|c|}
\hline
Task & Classifier & Transformation & \qupid-\fft{}  & \qupid-\db{1} & \qupid-\db{2} & \qupid-\db{3} & \qupid-\coif{1} & \qupid-\coif{2} & \qupid-\coif{3} & \qupid-\coif{4} & \qupid-\Id{} \\
\hline
\multirow{4}{*}{$\texttt{CD68}^+$ vs $\texttt{CD8}^+$} & \texttt{RFC} & PD & 69.1 & 71.6 & 74.1 & 76.5 & 74.6 & 74.8 & 77.3 & \textbf{79.0} & 69.9 \\
 & \texttt{LDA} & PD & 63.4 & 64.5 & 64.3 & 60.0 & 61.1 & \textbf{66.3} & 59.5 & 57.7 & 63.2 \\
 & \texttt{RFC} & Raw & 52.1 & 55.8 & 51.8 & 53.6 & 58.2 & 58.6 & 59.5 & \textbf{60.2} & 52.7 \\
 & \texttt{LDA} & Raw & 48.4 & 62.5 & 62.8 & 64.5 & 61.4 & \textbf{65.0} & 63.9 & 60.0 & 61.4 \\
\hline
\multirow{4}{*}{$\texttt{CD68}^+$ vs $\texttt{FoxP3}^+$} & \texttt{RFC} & PD & 94.3 & 91.7 & 94.3 & 94.8 & 96.2 & 96.8 & \textbf{98.2} & 97.0 & 91.7 \\
& \texttt{LDA} & PD & \textbf{94.5} & 87.8 & 86.5 & 90.3 & 88.0 & 91.5 & 88.7 & 89.0 & 90.0 \\
& \texttt{RFC} & Raw & 87.2 & 89.0 & 87.7 & 86.0 & 87.8 & 84.7 & 84.8 & 88.2 & \textbf{89.2} \\
 & \texttt{LDA} & Raw & 64.3 & 85.8 & 86.8 & 84.8 & 86.8 & 85.3 & \textbf{88.5} & 86.5 & 88.3 \\
\hline
\multirow{4}{*}{$\texttt{CD8}^+$ vs $\texttt{FoxP3}^+$} & \texttt{RFC} & PD & 91.7 & 89.1 & 89.1 & 90.5 & 94.5 & 92.7 & 95.2 & \textbf{96.3} & 90.7 \\
& \texttt{LDA} & PD & \textbf{85.2} & 80.9 & 80.7 & 77.7 & 80.7 & 82.3 & 80.0 & 80.7 & 82.7 \\
& \texttt{RFC} & Raw & 83.0 & 86.6 & 87.1 & 86.6 & 84.5 & 85.9 & 85.7 & \textbf{87.7} & 83.8 \\
 & \texttt{LDA} & Raw & 68.0 & 82.1 & 84.0 & 84.3 & 85.4 & 84.8 & 84.3 & \textbf{87.0} & 85.9 \\
\hline
\multirow{4}{*}{All classes} & \texttt{RFC} & PD & 72.6 & 67.4 & 72.7 & 74.7 & 78.1 & 78.0 & 77.5 & \textbf{79.4} & 68.8 \\
& \texttt{LDA} & PD & \textbf{67.9} & 59.5 & 61.7 & 58.4 & 61.4 & 60.3 & 56.1 & 56.9 & 47.1 \\
 & \texttt{RFC} & Raw & 59.2 & 60.8 & 60.2 & 62.9 & 62.7 & 62.7 & 63.8 & \textbf{66.2} & 62.0 \\
 & \texttt{LDA} & Raw & 41.0 & 64.9 & 67.0 & 65.6 & 65.6 & 66.5 & 67.7 & 67.2 & \textbf{67.9} \\
\hline
\end{tabular}}
\caption{Comparative analysis of immune cells classification for several choices of discrete transforms on Persistence diagrams and raw points with \texttt{RandomForestClassifier}  (\texttt{RFC}) and \texttt{LinearDiscriminantAnalysis} (\texttt{LDA}).}
\label{tab:complete:cancer}
\end{table*}
\end{document}